\documentclass{amsart}
\usepackage{url}
\usepackage{graphicx}
\newtheorem{thm}{Theorem}[section]

\newtheorem{lem}[thm]{Lemma}

\theoremstyle{definition}

\theoremstyle{remark}

\begin{document}
\title[Some characterizations of the quasi-sum production models]{Some characterizations of the quasi-sum production models with proportional marginal rate of substitution}
\author[A.D. V\^{\i}lcu, G.E. V\^{\i}lcu]{Alina Daniela V\^{\i}lcu, Gabriel Eduard V\^{\i}lcu}

\date{}
\maketitle

\abstract  In this note we classify quasi-sum production functions with constant elasticity of production with respect to any factor of production and with proportional marginal rate of substitution. 
\endabstract

The notion of \emph{production function} is a key concept in both macroeconomics and microeconomics, being used in the mathematical modeling of the relationship between the output of a firm, an industry, or an entire economy, and the inputs that have been used in obtaining it. Generally,  production function is a twice differentiable mapping $f:\mathbb{R}^n_+\rightarrow\mathbb{R}_+$, $f=f(x_1,\ldots,x_n)$,
where $f$ is the quantity of output, $n$ is the number of the inputs
and $x_1,\ldots,x_n$ are the factor inputs. A production function $f$ is called \emph{quasi-sum} \cite{AMih,CH4} if there are strict monotone functions $G,h_1,\ldots,h_n$ with $G'>0$ such that
\begin{equation}\label{1}
f(x)=G(h_1(x_1)+\ldots+h_n(x_n)),
\end{equation}
where $x=(x_1,\ldots,x_n)\in\mathbb{R}^n_+$. We note that these functions are of great interest because they appear as solutions of the general bisymmetry equation, being related to the problem of consistent aggregation \cite{AM}.

Among the family of production functions, the most famous is the so-called Cobb-Douglas production function. A generalized Cobb-Douglas production function depending on $n$-inputs is given by
\begin{equation}\label{2}
f(x_1,...,x_n)=A\cdot
\prod_{i=1}^n x_i^{\alpha_{i}},
\end{equation}
where $A,\alpha_1,\ldots,\alpha_n>0$.
We recall that a production function of the form
$
f(x)=G(h(x_1, . . . , x_n)),
$
where $G$ is a strictly increasing function and $h$ is a homogeneous function of any
given degree $p$, is said to be a \emph{homothetic} production function \cite{CH2}.
It is easy to see that a production function $f$
can be identified with the graph of $f$, \emph{i.e.} the nonparametric hypersurface of $\mathbb{E}^{n+1}$ defined by
\begin{equation}\label{3}
L(x_1,...,x_n)=(x_1,...,x_n,f(x_1,...,x_n))
\end{equation}
and called the \emph{production hypersurface} of $f$ (see \cite{VV,VGE}). Motivated by some recent classification results concerning production hypersurfaces \cite{AE3,CH4,CH2,FW,WF}, in the present work we classify quasi-sum production functions with proportional marginal rate of substitution and investigate the existence of such production models whose production hypersurfaces have null Gauss-Kronecker curvature or null mean curvature.
We recall that, if $f$ is a production function with $n$ inputs $x_1,x_2,...,x_n$, $n\geq 2$, then the \emph{elasticity of production} with respect to a certain factor of production $x_i$ is defined as
\begin{equation}\label{4}
E_{x_i}=\frac{x_i}{f}f_{x_i}
\end{equation}
and the \emph{marginal rate of technical substitution} of input $x_j$
for input $x_i$ is given by
\begin{equation}\label{5}
{\rm MRS}_{ij}=\frac{f_{x_j}}{f_{x_i}},
\end{equation}
where the subscripts denote partial derivatives of the function $f$
with respect to the corresponding variables.
A production function satisfies \emph{the proportional marginal rate of substitution property} if
\begin{equation}\label{6}
{\rm MRS}_{ij}=\frac{x_i}{x_j},\ {\rm for\ all\ } 1\leq i\neq j\leq n.
\end{equation}
In the last section of the paper we will prove the following theorem which generalize the results from \cite{VV2}.
\begin{thm}\label{T1}
Let $f$ be a quasi-sum production function given by {\rm (\ref{1})}.  Then:
\begin{enumerate}
  \item[i.] The elasticity of production is a constant $k_i$ with respect to a certain factor of production $x_i$ if and only if $f$ reduces to
\begin{equation}\label{7}
f(x_1,\ldots,x_n)=A\cdot x_i^{k_i}\cdot\exp\left(D\displaystyle\sum_{j\neq i}h_j(x_j)\right),
\end{equation}
where $A$ and $D$ are positive constants.
\item[ii.] The elasticity of production is a constant $k_i$ with respect to all factors of production $x_i$, $i=1,\ldots,n$,  if and only if $f$ reduces to the generalized Cobb-Douglas production function given by (\ref{2}).
  \item[iii.] The production function satisfies the proportional marginal rate of substitution property  if and only if it reduces to
  the homothetic generalized Cobb-Douglas production function given by
\begin{equation}\label{8}
f(x_1,\ldots,x_n)=F\left(\prod_{i=1}^n x_i^k\right),
\end{equation}
where $k$ is a nonzero real number.
\item[iv.] If the production function satisfies the proportional marginal rate of substitution property, then:
\begin{enumerate}
  \item[iv$_1$.] The production hypersurface has vanishing Gauss-Kronecker curvature if and only if, up to a suitable translation, $f$ reduces to the following generalized Cobb-Douglas production function with constant return to scale:
    \begin{equation}\label{9}
    f(x_1,\ldots,x_n)=A\cdot \prod_{i=1}^n x_i^{\frac{1}{n}}.
    \end{equation}
  \item[iv$_2$.] The production hypersurface cannot be minimal.
  \item[iv$_3$.] The production hypersurface has vanishing sectional curvature if and only if, up to a suitable translation, $f$ reduces to the following generalized Cobb-Douglas production function:
    \begin{equation}\label{10}
    f(x_1,\ldots,x_n)=A\cdot \prod_{i=1}^n \sqrt{x_i}.
    \end{equation}
\end{enumerate}
\end{enumerate}
\end{thm}



\section{Preliminaries on the geometry of hypersurfaces}
\label{PREL}

For general references on the geometry of hypersurfaces, we refer to \cite{CH1}.

If $M$ is a hypersurface of the Euclidean space $\mathbb E^{n+1}$, then it is known that the \emph{Gauss map}
$\nu: M \rightarrow S^{n}$ maps $M$ to the unit hypersphere
$S^n$ of $\mathbb E^{n+1}$. With the help of the differential $d\nu$ of $\nu$ it can be defined a linear operator on the tangent space $T_pM$, denoted by $S_p$ and known as the \emph{shape operator}, by
$
g(S_pv,w)=g(d\nu(v),w),
$
for $v,w\in T_pM$, where $g$ is the metric tensor on $M$ induced from the Euclidean metric
on $\mathbb E^{n+1}$.
The eigenvalues of the shape
operator are called\emph{ principal curvatures}. The determinant of the shape operator
$S_p$, denoted by $K(p)$, is called the \emph{Gauss-Kronecker curvature}. When $n = 2$, the
Gauss-Kronecker curvature is simply called the \emph{Gauss curvature}, which is
intrinsic due to famous Gauss's Theorema Egregium. The trace of
the shape operator $S_p$ is called the \emph{mean curvature} of the
hypersurfaces. In contrast to the Gauss-Kronecker curvature, the mean
curvature is extrinsic, which depends on the immersion of
the hypersurface. A hypersurface is said to be \emph{minimal }if its mean
curvature vanishes identically.
We recall now the following Lemma which will be used  in the proof of Theorem \ref{T1}.
\begin{lem}\label{L1} {\rm \cite{CH1}}
For the production hypersurface defined by {\rm (\ref{3})} and
$
w=\sqrt{1+\displaystyle\sum_{i=1}^{n} f_{i}^{2}},
$
we have:
\begin{enumerate}
  \item[i.] The Gauss-Kronecker curvature $K$ is given by
      \begin{equation}\label{11}
     K=\frac{\det(f_{x_ix_j})}{w^{n+2}}.
      \end{equation}
  \item[ii.] The mean curvature $H$ is given by
       \begin{equation}\label{11b}
     H=\frac{1}{n}\sum_{i=1}^{n}\frac{\partial}{\partial x_i}\left(\frac{f_{x_i}}{w}\right).
      \end{equation}
  \item[iii.] The sectional curvature $K_{ij}$ of the plane section spanned by $\frac{\partial}{\partial x_i}$, $\frac{\partial}{\partial x_j}$ is
        \begin{equation}\label{12}
      K_{ij}=\frac{f_{x_ix_i}f_{x_jx_j}-f^2_{x_ix_j}}{w^2\left(1+f_{x_i}^2+f_{x_j}^2\right)}.
      \end{equation}
\end{enumerate}
\end{lem}

\section{Proof of Theorem \ref{T1}}

Let $f$ be a quasi-sum production function given by (\ref{1}). Then we have
 \begin{equation}\label{14}
f_{x_i}(x)=G'(u) h'_i(x_i)
\end{equation}
with $u=h_1(x_1)+\ldots+h_n(x_n)$ and from (\ref{14}) we derive
 \begin{equation}\label{15}
f_{x_ix_i}=G''(h'_i)^2+G'h''_i,\ i=1,\ldots,n,
\end{equation}
 \begin{equation}\label{16}
f_{x_ix_j}=G''h'_ih'_j,\ i\neq j.
\end{equation}

i. We first prove the left-to-right implication. If the elasticity of production is a constant $k_i$ with respect to a certain factor of production $x_i$, then from (\ref{4}) we obtain
 \begin{equation}\label{17}
f_{x_i}=k_i\frac{f}{x_i}.
\end{equation}

Using now (\ref{1}) and (\ref{14}) in (\ref{17}) we get
 \begin{equation}\label{18}
\frac{G'}{G}=k_i\frac{1}{x_ih'_i}.
\end{equation}

By taking the partial derivative of (\ref{18}) with respect to $x_j$, $j\neq i$, we obtain
\[
h'_j\frac{G''G-(G')^2}{G^2}=0.
\]

Now, taking into account that $h_j$ is a strict monotone function, we find
\begin{equation}\label{20}
G(u)=C\cdot e^{Du},
\end{equation}
for some positive constants $C$ and $D$.
Hence from (\ref{18}) and (\ref{20}) we obtain
\begin{equation}\label{21}
h_i(x_i)=\frac{k_i}{D}\ln x_i+A_i,
\end{equation}
where $A_i$ is a real constant. Finally, combining (\ref{1}), (\ref{20}) and (\ref{21}) we get a function of the form (\ref{7}), where $A=Ce^{D\cdot A_i}$.
The converse can be verified easily by direct computation.

ii. The assertion is an immediate consequence of i.

iii. Assume first that $f$ satisfies the proportional marginal rate of substitution property. Then from (\ref{5}), (\ref{6}) and (\ref{14}) we derive
$
x_ih'_i=x_jh'_j,\ \forall i\neq j.
$
Hence we conclude that there exists a nonzero real number $k$ such that:
$
x_ih'_i=k,\ i=1,\ldots,n,
$
and therefore we obtain
\begin{equation}\label{22}
h_i(x_i)=k\ln x_i+C_i,\ i=1,\ldots,n,
\end{equation}
for some real constants $C_1,\ldots,C_n$. Now,
from (\ref{1}) and (\ref{22}) we derive
\[
f(x)=G\left(k\displaystyle\sum_{i=1}^{n}\ln x_i+\overline{A}\right),
\]
where $\overline{A}=\displaystyle\sum_{i=1}^{n}C_i$ and hence we find
\begin{equation}\label{23}
f(x)=(G\circ\ln) \left(A\cdot\prod_{i=1}^n x_i^k\right),
\end{equation}
where $A=e^{\overline{A}}$. Therefore we get a production function of the form (\ref{8}), where $F(u)=(G\circ\ln)(A\cdot u)$.

The converse is easy to verify.

iv$_1$. We first prove the left-to-right implication. If the production hypersurface has null Gauss-Kronecker curvature, then from (\ref{11}) we get
\begin{equation}\label{23a}
\det(f_{x_ix_j})=0.
\end{equation}

On the other hand, the determinant of the Hessian matrix of $f$ is given by \cite{CH8}
\begin{equation}\label{23b}
\det(f_{x_ix_j})=(G')^n\prod_{i=1}^n h_i''+(G')^{n-1}G''\sum_{i=1}^nh_1''\cdot\ldots\cdot h_{i-1}''(h_i')^2h_{i+1}''\cdot\ldots\cdot h_n''.
\end{equation}

By using (\ref{22}), (\ref{23a}) and (\ref{23b}), we obtain
\[
(-1)^n(G')^{n-1}k^n(G'-knG'')=0.
\]

But $G'>0$ and $k\neq 0$ and hence we derive
\begin{equation}\label{23c}
\frac{G''}{G'}=\frac{1}{kn}.
\end{equation}

After solving (\ref{23c}) we find
\begin{equation}\label{23d}
G(u)=Cnke^{\frac{u}{nk}}+D
\end{equation}
for some constants $C,D$ with $C>0$. Combining (\ref{23}) and (\ref{23d}), after a suitable translation, we conclude that the function $f$ reduces to the
form (\ref{9}).
The converse follows easily by direct computation.

iv$_2$. Let us assume that the production hypersurface is minimal. Then we have $H=0$ and from  (\ref{11b}) we derive
\[
\sum_{i=1}^{n}f_{x_ix_i}\left(1+\sum_{i=1}^{n}f_{x_i}^2\right)-\sum_{i,j=1}^{n}f_{x_i}f_{x_j}f_{x_ix_j}=0
\]
which reduces to
\begin{equation}\label{25}
\sum_{i=1}^{n}f_{x_ix_i}+\sum_{i\neq j}\left(f_{x_i}^2f_{x_jx_j}-f_{x_i}f_{x_j}f_{x_ix_j}\right)=0.
\end{equation}

By introducing (\ref{14}), (\ref{15}) and (\ref{16}) in (\ref{25}), we get
\begin{equation}\label{26}
G''\sum_{i=1}^{n}(h'_i)^2+G'\sum_{i=1}^{n}h''_i+(G')^3\sum_{i\neq j}(h'_i)^2h''_j=0.
\end{equation}

By using now (\ref{22}) in (\ref{26}) and taking into account that $k\neq 0$, we obtain
\begin{equation}\label{27}
(kG''-G')\sum_{i=1}^{n}\frac{1}{x_i^2}-k^2(G')^3\sum_{i\neq j}\frac{1}{x_i^2x_j^2}=0.
\end{equation}

But the only solution of the equation (\ref{27}) is $G(u)=constant$, which is a contradiction because $G'>0$. Hence the production hypersurface cannot be minimal.

iv$_3$. Assume first that the production hypersurface has $K_{ij}=0$. Then from (\ref{12}) we get
\begin{equation}\label{32}
f_{x_ix_i}f_{x_jx_j}-f^2_{x_ix_j}=0.
\end{equation}

By introducing (\ref{14}), (\ref{15}) and (\ref{16}) in (\ref{32}), since $G'\neq 0$, we obtain
\begin{equation}\label{33}
[(h'_i)^2h_j''+(h'_j)^2h''_i]G''+h''_ih''_jG'=0.
\end{equation}

By using now (\ref{22}) in (\ref{33}) and taking into account that $k\neq 0$, we obtain
\begin{equation}\label{34}
\frac{G''}{G'}=\frac{1}{2k}.
\end{equation}

After solving (\ref{34}) we get
\begin{equation}\label{35}
G(u)=2kCe^{\frac{u}{2k}}+D
\end{equation}
for some constants $C,D$ with $C>0$.
Finally, combining (\ref{23}) and (\ref{35}), after a suitable translation, we conclude that the function $f$ reduces to the
Cobb-Douglas production function given by (\ref{10}).
The converse is easy to verify by direct computation.

\section*{Acknowledgements}
The second author was supported by National Research Council - Executive Agency for Higher Education Research and Innovation Funding
(CNCS-UEFISCDI), project number PN-II-ID-PCE-2011-3-0118.

Alina Daniela V\^{I}LCU\\
      Petroleum-Gas University of Ploie\c sti,\\
      Department of Computer Science, Information Technology, Mathematics and Physics,\\
      Bulevardul Bucure\c sti, Nr. 39, Ploie\c sti 100680-ROMANIA\\
      e-mail: daniela.vilcu@upg-ploiesti.ro\\

Gabriel Eduard V\^{I}LCU$^{1,2}$ \\
      $^1$University of Bucharest, Faculty of Mathematics and Computer Science,\\
      Research Center in Geometry, Topology and Algebra,\\
      Str. Academiei, Nr. 14, Sector 1,\\
      Bucure\c sti 70109-ROMANIA\\
      e-mail: gvilcu@gta.math.unibuc.ro\\
      $^2$Petroleum-Gas University of Ploie\c sti,\\
       Department of Mathematical Modelling, Economic Analysis and Statistics,\\
      Bulevardul Bucure\c sti, Nr. 39, Ploie\c sti 100680-ROMANIA\\
      e-mail: gvilcu@upg-ploiesti.ro

\end{document}